\documentclass[10pt,a4paper]{amsart}
\date{29 January 2014}      
\usepackage{hyperref}
\usepackage{latexsym,amsmath,amsfonts,amscd,amssymb,graphics}
\usepackage[all]{xy}
\usepackage[ansinew]{inputenc}
\usepackage{float}
\usepackage{setspace}

\theoremstyle{plain}  
\newtheorem{theorem}{Theorem}[section]
\newtheorem*{theorem*}{Theorem}

\newtheorem{corollary}[theorem]{Corollary}
\newtheorem{lemma}[theorem]{Lemma}
\newtheorem{proposition}[theorem]{Proposition}

\theoremstyle{definition}
\newtheorem{definition}[theorem]{Definition}
\newtheorem{notation}[theorem]{Notation}
\newtheorem{assumption}[theorem]{Assumption}

\theoremstyle{remark}

\newtheorem{remark}[theorem]{Remark}

\newtheorem*{claim*}{Claim}

\numberwithin{equation}{section}

\newcommand{\suchthat}{\;|\;}
\newcommand{\st}{\;|\;}

\renewcommand{\leq}{\leqslant}

\renewcommand{\geq}{\geqslant}

\renewcommand{\setminus}{\smallsetminus}

\newcommand{\R}{\mathbb{R}}

\newcommand{\Z}{\mathbb{Z}}
\newcommand{\C}{\mathbb{C}}

\newcommand{\PP}{\mathbb{P}}

\newcommand{\cM}{\mathcal{M}}
\newcommand{\cN}{\mathcal{N}}

\newcommand{\cS}{\mathcal{S}}

\newcommand{\Sp}{\mathrm{Sp}}

\DeclareMathOperator{\Jac}{Jac}

\DeclareMathOperator{\divisor}{div}

\DeclareMathOperator{\cSym}{Sym}

\DeclareMathOperator{\rk}{rk}

\DeclareMathOperator{\End}{End}

\DeclareMathOperator{\codim}{codim}

\newcommand{\Pic}{\operatorname{Pic}}

\newcommand{\Aut}{\operatorname{Aut}}

\usepackage[top=25mm,bottom=25mm,left=25mm,right=25mm]{geometry}

\let\oldmarginpar\marginpar
\renewcommand\marginpar[1]{\oldmarginpar{\tiny\bf\begin{flushleft} #1
\end{flushleft}}}

\begin{document}

\title[Torelli theorem for the moduli spaces of rank 2 quadratic pairs]{Torelli theorem for the moduli spaces  of rank 2 quadratic pairs}

\author[A. G. Oliveira]{Andr\'e G. Oliveira}
\address{Departamento de Matem\'atica\\
  Universidade de Tr\'as-os-Montes e Alto Douro (UTAD)\\
  Quinta dos Prados \\ 5000-911 Vila Real \\ Portugal\\ www.utad.pt }
\email{agoliv@utad.pt}

\thanks{
The author thanks I. Biswas, P. Gothen and V. Muñoz for enlightening discussions. Also thanks the Isaac Newton Institute in Cambridge and Centre de Recerca Matematica in Barcelona --- that he visited while preparing the paper --- for the excellent conditions provided.
Member of VBAC (Vector Bundles on Algebraic Curves) and GEAR (Geometric Structures and Representation Varieties)
Partially supported the FCT (Portugal) with national funds 
through the projects PTDC/MAT/099275/2008, PTDC/MAT/098770/2008, PTDC/MAT/120411/2010 and PTDC/MAT-GEO/0675/2012
and through Centro de Matem\'atica da
Universidade de Tr\'as-os-Montes e Alto Douro (PEst-OE/MAT/UI4080/2011). 
}


\begin{abstract}
Let $X$ be a smooth projective complex curve. We prove that a Torelli type theorem holds, under certain conditions, for the moduli space of $\alpha$-polystable quadratic pairs on $X$ of rank 2.
\end{abstract}

\maketitle


\section{Introduction}

Let $X$ be a smooth projective curve over $\C$. This paper deals with moduli spaces of quadratic pairs $(V,\gamma)$ consisting of a vector bundle $V$ together with a non-zero bilinear symmetric map $\gamma$ on $V$, with values in a fixed line bundle $U$.
The semistability condition for these objects depends on a real parameter $\alpha$.
Let $\cN_\alpha^\Lambda$ denote the moduli space of $\alpha$-polystable $U$-quadratic pairs $(V,\gamma)$ with $\rk(V)=2$ and $\deg(V)=d$, where the determinant of $V$ is fixed to be $\Lambda\in\Jac^d(X)$.
The goal of this paper is to give new insights on the geometry and topology of $\cN_\alpha^\Lambda$. The main result is a proof that, under some conditions, a Torelli type theorem holds for $\cN_\alpha^\Lambda$  --- see Theorem \ref{torelli for quad pairs} for the precise statement.
We follow the methods used by Muñoz in \cite{munoz:2009}, with the help of the results of \cite{gothen-oliveira:2011}. 
 One of the differences is that, in our case, the moduli spaces are not smooth. 

The fact that $\gamma$ is a non-linear object (as opposed to most studied cases of vector bundles with extra structure) is perhaps one of the reasons why these moduli spaces have not yet been much studied. In fact, this non-linearity of $\gamma$ is in the origin of many difficulties which arise, starting with the general definition of stability. For this reason, we consider only the rank $2$ case.
If $d=d_U$, where $d_U$ denotes the degree of $U$, the dependence of the semistability condition on $\alpha$ disappears, and we get the usual stability condition for orthogonal bundles. In other words, we get the moduli space of rank $2$ orthogonal vector bundles. These moduli spaces constitute hence natural generalisations of moduli spaces of orthogonal bundles, so it is natural to study their geometry even in this low rank case. 
If $d_U>d$ the moduli spaces are non-empty for $\alpha\leq d/2$, and their dimension increases with the difference $d_U-d$. For example, it is known \cite{serman:2008} that the forgetful map $(V,\gamma)\mapsto V$ is an embedding of the moduli space of orthogonal bundles into the moduli space vector bundles; however, when the difference $d_U-d$ is at least equal to the genus of $X$, then the same forgetful map is surjective and generically a projective bundle. This forgetful map will play a central role in this paper.

These spaces also arise naturally as special subvarieties of the moduli spaces of $\mathrm{Sp}(2n,\mathbb{R})$-Higgs bundles \cite{garcia-prada-mundet:2004,gothen-oliveira:2011}, so the study of its geometry and topology is useful to draw conclusions on the topology of character varieties for representations of surface groups in $\mathrm{Sp}(2n,\mathbb{R})$.

\section{Quadratic pairs and their moduli spaces}
\subsection{Quadratic pairs}

Let $X$ be a smooth projective curve over $\C$ of genus $g\geq 2$. Let $U$ be a once and for all fixed holomorphic line bundle over
$X$, and let $d_U=\deg(U)$ denote the degree of $U$.

\begin{definition}\label{definition of U quadratic pair}
  A \emph{$U$-quadratic pair} on $X$ is a pair $(V,\gamma)$, where $V$
  is a holomorphic vector bundle over $X$ and $\gamma$ is a global
  holomorphic non-zero section of $S^2V^*\otimes U$, i.e., $\gamma\in
  H^0(X,S^2V^*\otimes U)$.
\end{definition}

Of course $\gamma$ can also be seen as a holomorphic map $\gamma:V\to V^*\otimes U$ which is symmetric, i.e.,
$\gamma^t\otimes 1_U=\gamma$.
 The \emph{rank} and \emph{degree} of a quadratic pair $(V,\gamma)$ are the rank and degree of $V$. We shall say that $(V,\gamma)$ is of \emph{type}
  $(n,d)$ if $\rk(V)=n$ and $\deg(V)=d$.
We will often refer to a $U$-quadratic pair simply as a
\emph{quadratic pair}. 
Two $U$-quadratic pairs $(V,\gamma)$ and $(V',\gamma')$ are
  \emph{isomorphic} if there exists an isomorphism $f:V\to V'$ such that
  $\gamma'f=((f^t)^{-1}\otimes 1_U)\gamma$.
For instance, for any $\lambda\in\C^*$, the pairs $(V,\gamma)$ and $(V,\lambda\gamma)$ are isomorphic.

\subsection{\texorpdfstring{$U$}{U}-quadratic pairs of rank \texorpdfstring{$1$}{1}}

Our main object of study will be the moduli space $\alpha$-polystable quadratic pairs of type $(2,d)$, with fixed determinant. However, we will occasionally have to consider the moduli space of quadratic pairs of type $(1,d)$. 
These kind of pairs, being given by a line bundle $L$ together with a section of the line bundle $L^{-2}U$, are quite simple objects.
  Fix a real parameter $\alpha$. A $U$-quadratic pair $(L,\delta)$ of
  type $(1,d)$ is \emph{$\alpha$-stable} if $\alpha\leq d$.
Denote the moduli space of $\alpha$-stable $U$-quadratic pairs on $X$ of rank $n$ and degree $d$ by $\cN_\alpha(1,d)$.
For a fixed type $(1,d)$, all the moduli spaces
$\cN_\alpha(1,d)$ with $\alpha\leq d'$ are isomorphic and there is
only one so-called \emph{critical value} of $\alpha$, for which the
stability condition changes, namely $\alpha=d$.
$\cN_\alpha(1,d)$ is a particularly simple space.
The following proposition summarizes its features. For details, see Section 2.2 of \cite{gothen-oliveira:2011}.

\begin{proposition}\label{prop:Nalpha(1,d)}
Let $\cN_\alpha(1,d)$ be the moduli space of $\alpha$-stable
quadratic pairs of type $(1,d)$.
\begin{enumerate}
\item If $d>d_U/2$ or $\alpha>d$, then $\cN_\alpha(1,d)=\emptyset$;
\item If $d<d_U/2$ and $\alpha\leq d$, then $\cN_\alpha(1,d)$ is the
  $2^{2g}$-fold cover of the symmetric product $\cSym^{d_U-2d}(X)$
  obtained by pulling back, via the Abel-Jacobi map, the cover
  $\Jac^{d}(X)\to\Jac^{d_U-2d}(X)$ given $L\mapsto L^{-2}U$. In particular, $\cN_\alpha(1,d)$ is a smooth projective variety of dimension $d_U-2d$.
\item If $d_U$ is even and $\alpha\leq d_U/2$, then $\cN_\alpha(1,d_U/2)$ is the set of the $2^{2g}$ square roots of $U$.  
\end{enumerate}
\end{proposition}

\subsection{Stability of quadratic pairs of rank $2$}

Let us now consider $U$-quadratic pairs $(V,\gamma)$ of type $(2,d)$.
Given a rank $2$ vector bundle $V$ and a line subbundle $L\subset V$,
 denote by $L^\perp$ the kernel of the projection $V^*\to L^{-1}$. 
It is thus a line subbundle of $V^*$ and $V/L$ is isomorphic to
$(L^\perp)^{-1}$.

\begin{definition}\label{simplified notion of polystability}
Let $(V,\gamma)$ be a quadratic pair of type $(2,d)$.
\begin{itemize}
\item The pair $(V,\gamma)$ is \emph{$\alpha$-semistable} if
  $\alpha\leq d/2$ and, for any line bundle $L\subset V$, the
  following conditions hold:
\begin{enumerate}
\item $\deg(L)\leq\alpha$, if $\gamma(L)=0$;
\item $\deg(L)\leq d/2$, if $\gamma(L)\subset L^\perp U$;
\item $\deg(L)\leq d-\alpha$, if $\gamma(L)\not\subset L^\perp U$.
\end{enumerate}
\item The pair $(V,\gamma)$ is \emph{$\alpha$-stable} if it is
  $\alpha$-semistable for any line bundle $L\subset V$, the conditions
  $(1)$, $(2)$ and $(3)$ above hold with strict inequalities.
\item The pair $(V,\gamma)$ is \emph{$\alpha$-polystable} if it is $\alpha$-semistable and, for any line bundle $L\subset V$, the following conditions hold:
\begin{enumerate}
\item if $\gamma(L)=0$ and $\deg(L)=\alpha$, then there is an $L'\subset V$ such that $V=L\oplus L'$
  and with respect to this decomposition,
$$\gamma=\begin{pmatrix}
          0 & 0\\
	  0 & \gamma'
         \end{pmatrix}$$ with $\gamma'\in H^0(X,L'^{-2}U)$ non-zero;
\item if $\gamma(L)\subset L^\perp U$ and $\deg(L)=d/2$, then there is $L'\subset V$ such that $V=L\oplus L'$ and with respect to this decomposition,
$$\gamma=\begin{pmatrix}
          0 & \gamma'\\
	  \gamma' & 0
         \end{pmatrix}$$ with $\gamma'\in H^0(X,L^{-1}L'^{-1}U)$ non-zero;
\item if $\gamma(L)\not\subset L^\perp U$ and $\deg(L)=d-\alpha$, then there is $L'\subset V$ such that $V=L\oplus L'$ and with respect to this decomposition, 
$$\gamma=\begin{pmatrix}
          \gamma' & 0\\
	  0 & 0
         \end{pmatrix}$$ with $\gamma'\in H^0(X,L^{-2}U)$ non-zero.
\end{enumerate}
\end{itemize}
\end{definition}

It will be useful for us to classify the possible types of destabilizing subbundles:
\begin{definition}\label{def:types}
  Let $(V,\gamma)$ be a quadratic pair of type $(2,d)$, and let $L\subset V$ be an $\alpha$-destabilizing line subbundle of $(V,\gamma)$. We say that $L$ is:
  \begin{enumerate}
  \item of \emph{type \textbf{(A)}}, if $\deg(L)\geq \alpha$ and $\gamma(L)=0$.
  \item of \emph{type \textbf{(B)}}, if $\deg(L)\geq d/2$ and $\gamma(L)\subset L^{\perp}U$.
  \item of \emph{type \textbf{(C)}}, if $\deg(L)\geq d-\alpha$ and $\gamma(L)\not\subset
  L^{\perp}U$.
  \end{enumerate}
\end{definition}

%

For a given $\alpha\in\R$, denote by 
$\cN_\alpha=\cN_{X,\alpha}(2,d)$ the moduli space of $\alpha$-polystable $U$-quadratic pairs of type $(2,d)$ on the curve $X$.
We obtain the moduli space of $\alpha$-polystable $U$-quadratic pairs with fixed determinant in the usual way. There is a map 
\begin{equation}\label{eq:detmap}
\det:\cN_\alpha(2,d)\to\Jac^d(X), \hspace{.5cm}\det(V,\gamma)=\bigwedge\nolimits^2 V
\end{equation} and the moduli space of $\alpha$-polystable $U$-quadratic pairs with fixed determinant $\Lambda\in\Jac^d(X)$ is the subvariety of $\cN_\alpha(2,d)$ defined as $$\cN^\Lambda_\alpha=\cN^\Lambda_{X,\alpha}=\mathrm{det} ^{-1}(\Lambda).$$

From \cite[Theorem 1.6]{gomez-sols:2000}, $\cN_\alpha$ is a projective variety, and from \cite[Theorem 5.1]{gothen-oliveira:2011}, we see that, for $\alpha\leq d/2$, and $d<d_U$
\begin{equation}\label{eq:Nalpha}
\dim(\cN^\Lambda_\alpha)=3(d_U-d)-1
\end{equation}


\subsection{Deformation theory of quadratic pairs}
\label{sec:deformation-theory}

Here we look at the local structure of the moduli space of quadratic pairs. The deformation theory of a quadratic pair $(V,\gamma)\in\cN_\alpha$ is governed by the following complex of sheaves on $X$ (cf.  \cite{biswas-ramanan:1994}):
$$C^\bullet(V,\gamma):\End(V)\xrightarrow{\rho(\gamma)}S^2V^*\otimes U,\hspace{0,5cm}\rho(\gamma)(\psi)=-(\psi^t\otimes 1_U)\gamma-\gamma\psi.$$
Moreover, $\rho(\gamma)$ induces a long exact
sequence
\begin{equation}\label{eq:les}
\begin{split}
0&\longrightarrow\mathbb{H}^0(X,C^\bullet(V,\gamma))\longrightarrow
H^0(X,\End(V))\longrightarrow
H^0(X,S^2V^*\otimes U)\longrightarrow\\
&\longrightarrow\mathbb{H}^1(X,C^\bullet(V,\gamma))\longrightarrow H^1(X,\End(V))\longrightarrow
H^1(X,S^2V^*\otimes U)\longrightarrow\\
&\longrightarrow\mathbb{H}^2(X,C^\bullet(V,\gamma))\longrightarrow 0.
\end{split}
\end{equation}
So, $\mathbb{H}^0(X,C^\bullet(V,\gamma))$ can be identified with the infinitesimal automorphism space of $(V,\gamma)$:
$$\mathbb{H}^0(X,C^\bullet(V,\gamma))\cong\{\psi\in\End(V)\,\st\,\rho(\gamma)(\psi)=0\}.$$ Also, from \cite{biswas-ramanan:1994}, it is known that $\mathbb{H}^1(X,C^\bullet(V,\gamma))$ is canonically isomorphic to the space of infinitesimal deformations of $(V,\gamma)$. In particular, if $(V,\gamma)$ represents a smooth point of the moduli space $\cN_\alpha$, then $\mathbb{H}^1(X,C^\bullet(V,\gamma))$ is naturally isomorphic to the tangent space of $\cN_\alpha$ at the point represented by $(V,\gamma)$.

If we now restrict the deformation theory to the fixed determinant moduli space, and if $(V,\gamma)\in\cN^\Lambda_\alpha$, we have to consider the subcomplex $C^\bullet_0(V,\gamma)$ of $C^\bullet(V,\gamma)$ given by
$C^\bullet_0(V,\gamma):\End_0(V)\xrightarrow{\rho(\gamma)}S^2V^*\otimes U$, where $\End_0(V)$ denotes the sheaf of traceless endomorphisms of $V$.
The trace yields a short exact sequence 
$$0\longrightarrow\End_0(V)\longrightarrow\End(V)\longrightarrow\mathcal{O}\longrightarrow0$$ which in turn gives rise to an exact sequence of complexes, inducing the following long exact sequence: 
\begin{equation}\label{eq:sesfixeddet}
\begin{split}
0&\longrightarrow\mathbb{H}^0(X,C_0^\bullet(V,\gamma))\longrightarrow\mathbb{H}^0(X,C^\bullet(V,\gamma))\longrightarrow H^0(X,\mathcal{O})\longrightarrow\\
&\longrightarrow\mathbb{H}^1(X,C_0^\bullet(V,\gamma))\longrightarrow\mathbb{H}^1(X,C^\bullet(V,\gamma))\longrightarrow H^1(X,\mathcal{O})\longrightarrow\\
&\longrightarrow\mathbb{H}^2(X,C_0^\bullet(V,\gamma))\longrightarrow\mathbb{H}^2(X,C^\bullet(V,\gamma))\longrightarrow 0.
\end{split}
\end{equation}

\begin{remark}\label{rem:infdefspacefixeddet}
The infinitesimal deformation space of $(V,\gamma)$ in $\cN^\Lambda_\alpha$ is the kernel of the map $\mathbb{H}^1(X,C^\bullet(V,\gamma))\to H^1(X,\mathcal{O})$ in the above sequence. In particular, if $(V,\gamma)$ represents a smooth point of $\cN^\Lambda_\alpha$, then this kernel is the tangent space of $\cN^\Lambda_\alpha$ at $(V,\gamma)$.
\end{remark}

\begin{definition}\label{def of simple in quadratic pairs}
    A quadratic pair $(V,\gamma)$ is \emph{simple} if the group
  $\Aut(V,\gamma)$ of automorphisms of $(V,\gamma)$ is equal to $\{\pm
  1_V\}\cong\Z/2$.
\end{definition}


\begin{definition}\label{df:smooth}
Let $\cN^{\Lambda,sm}_\alpha$ be the subspace of $\cN^\Lambda_\alpha$ consisting of pairs which are $\alpha$-stable, simple and such that $\mathbb{H}^2(X,C_0^\bullet(V,\gamma))$ vanishes.
\end{definition}

By \cite[Proposition 2.21]{gothen-oliveira:2011}, $\cN^{\Lambda,sm}_\alpha$ is then a subspace of the smooth locus of $\cN^\Lambda_\alpha$, so from \eqref{eq:Nalpha} we have that, if $d<d_U$ and $\alpha\leq d/2$, then
\begin{equation}\label{eq:dim of smooth moduli of fixed det quad pairs}
\dim\cN^{\Lambda,sm}_\alpha=3(d_U-d)-1.
\end{equation}
Of course this follows as well from Remark \ref{rem:infdefspacefixeddet}.

\subsection{Changes with the parameter}
\label{sec:variation-of-moduli}

Here we give an overview about the variation of the moduli spaces $\mathcal{N}^\Lambda_{\alpha}$ with the stability parameter $\alpha$. This will be mainly an adaptation of the results obtained in \cite{gothen-oliveira:2011} to the fixed determinant moduli space.

It is clear that the semistability condition only changes at a finite number of values of the parameter $\alpha$. These values of $\alpha$ are called \emph{critical values}, and if $\alpha$ is not critical, we say that it is \emph{generic}. So, if we only allow to change $\alpha$ between any two consecutive critical values, all the moduli spaces $\cN^\Lambda_\alpha$ will be isomorphic.

Given a real number $x$, denote by $\lfloor x\rfloor$ the greatest integer smaller or equal than $x$.
The critical values of the parameter $\alpha$ are (see \cite[Proposition 3.3]{gothen-oliveira:2011}) all the integers between $d-\lfloor d_U/2\rfloor$ and $\lfloor d/2\rfloor$, together with $d/2$.

\begin{notation}\label{formula alphak}
  For each integer $d-\lfloor d/2\rfloor-\lfloor d_U/2\rfloor\leq k\leq 0$, define
  $\alpha_k=\lfloor d/2\rfloor+k$. Also, let 
  $\alpha_M=d/2$ and $\alpha_m=\alpha_{d-\lfloor d/2\rfloor-\lfloor d_U/2\rfloor}=d-\lfloor d_U/2\rfloor$
  so that $\alpha_M$ is the biggest critical value and $\alpha_m$ is the smallest one. Finally, let $\alpha_k^+$ denote the value of any parameter between the
  critical values $\alpha_k$ and $\alpha_{k+1}$, and let $\alpha_k^-$
  denote the value of any parameter between the critical values
  $\alpha_{k-1}$ and $\alpha_k$.
\end{notation}

With this notation, we have $\cN^\Lambda_{\alpha_k^+}=\cN^\Lambda_{\alpha_{k+1}^-}$
for all critical values $\alpha_k$.

%

\begin{definition}\label{def:flip loci}
  For each $k\in\left\{d-\lfloor d/2\rfloor-\lfloor d_U/2\rfloor,\ldots,0\right\}$, let
  $\cS^\Lambda_{\alpha_k^+}$ be the locus in $\cN^\Lambda_{\alpha_k^+}$ of \emph{pairs which are
  $\alpha_k^+$-semistable but $\alpha_k^-$-unstable}, that is,
   $\cS^\Lambda_{\alpha_k^+}=\left\{(V,\gamma)\in\cN^\Lambda_{\alpha_k^+}
    \suchthat(V,\gamma)\notin\cN^\Lambda_{\alpha_k^-}\right\}.$
  Similarly, define $\cS^\Lambda_{\alpha_k^-}$ to be the locus in $\cN^\Lambda_{\alpha_k^-}$ of \emph{pairs which are $\alpha_k^-$-semistable but
  $\alpha_k^+$-unstable},
  $\cS^\Lambda_{\alpha_k^-}=\left\{(V,\gamma)\in\cN^\Lambda_{\alpha_k^-}
    \suchthat(V,\gamma)\notin\cN^\Lambda_{\alpha_k^+}\right\}.$
  The spaces $\cS^\Lambda_{\alpha_k^\pm}$ are called the \emph{flip loci}
  for the critical value $\alpha_k$.
\end{definition}

The same arguments as the ones used in Propositions 3.6 and 3.7 and in Remark 3.8 of \cite{gothen-oliveira:2011} apply also to $\cS^\Lambda_{\alpha_k^\pm}$, showing that $\cS^\Lambda_{\alpha_k^\pm}$ lies in the $\alpha_k^\pm$-stable locus of $\cN^\Lambda_{\alpha_k^\pm}$ and that it is indeed a subvariety of $\cN^\Lambda_{\alpha_k^\pm}$.
An immediate consequence of the definition of the flip loci is that
$\cN^\Lambda_{\alpha_k^+}-\cS^\Lambda_{\alpha_k^+}
  \cong\cN^\Lambda_{\alpha_k^-}-\cS^\Lambda_{\alpha_k^-}$.

\vspace{.5cm}

Let us now give a brief description of the geometry of the flip loci.
Of course $\cS^\Lambda_{\alpha_M^+}=\emptyset$. With respect to the flip loci
$\cS^\Lambda_{\alpha_k^+}$ for the other critical values, we have the following, whose proof is entirely analogous to the one of \cite[Proposition 3.10]{gothen-oliveira:2011}.

\begin{proposition}\label{S+}
Let $\alpha_k<d/2$. Then $\cS^\Lambda_{\alpha_k^+}$ is a projective bundle over $\cN_{d-\alpha_k^+}(1,d-\alpha_k)$, with fiber $\mathbb{P}^{d-2\alpha_k+g-2}$.
\end{proposition}

From the previous proposition and from (2) of Proposition \ref{prop:Nalpha(1,d)} (in the case of $\cS^\Lambda_{\alpha_m^+}$ and $d_U$ even, use instead (3) of the same proposition), we have that, for every $\alpha_k<d/2$,
\begin{equation}\label{Salphak+(2,Lambda)}
\dim\cS^\Lambda_{\alpha_k^+}=d_U-d+g-2.
\end{equation}

Now we turn our attention to the other flip loci, $\cS^\Lambda_{\alpha_k^-}$. As in the case of $\cS^\Lambda_{\alpha_k^+}$, the behaviour of $\cS^\Lambda_{\alpha_k^-}$ depends on whether $\alpha_k=\alpha_M=d/2$ or not.
Since $\cN^\Lambda_\alpha=\emptyset$ for $\alpha>\alpha_M$, then 
\begin{equation}\label{eq:SalphaM-}
\cS^\Lambda_{\alpha_M^-}=\cN_{\alpha_M^-},
\end{equation} and the next result gives a description of this space in terms of the moduli space of  vector bundles.
Let $\cM^\Lambda=\cM_X(2,\Lambda)$ be the moduli space of rank $2$ polystable vector bundles, over $X$, of degree $d$ and fixed determinant $\Lambda\in\Jac^d(X)$.

\begin{assumption}\label{assump:dU-d>2g-2}
From now on, we assume that $d_U-d>2g-2$. This fact will be systematically used throughout the paper.
\end{assumption}

The next result follows from is the fixed determinant version of \cite[Proposition 3.13]{gothen-oliveira:2011}.
\begin{proposition}\label{fiber d/2}
Let $d$ and $d_U$ be as in Assumption \ref{assump:dU-d>2g-2}. Then the forgetful map 
$\pi:\cN_{\alpha_M^-}\to\cM^\Lambda$, $\pi(V,\gamma)=V$, is surjective and the fiber over a stable vector bundle $V$ is $\mathbb P^{3(d_U-d-g)+2}$.
\end{proposition}

We now give the description of the flip loci $\cS^\Lambda_{\alpha_k^-}$ with $\alpha_k<\alpha_M=d/2$. These are harder to describe than $\cS^\Lambda_{\alpha_k^+}$, and here we will only briefly sketch the description of $\cS^\Lambda_{\alpha_k^-}$; all the details can be found in \cite{gothen-oliveira:2011}.
If $(V,\gamma)\in\cS^\Lambda_{\alpha_k^-}$, then it is $\alpha_k^-$-stable and $\alpha_k^+$-unstable hence strictly $\alpha_k$-semistable. The destabilizing subbundle must be a line subbundle $M\subset V$ such that 
\begin{equation}\label{eq:MdestabS-}
\deg(M)=d-\alpha_k\hspace{.5cm}\text{ and }\hspace{.5cm}\gamma(M)\not\subset M^\perp U
\end{equation} hence $\gamma$ induces a non-zero holomorphic section
\begin{equation}\label{gamma'}
\gamma'\in H^0(X,M^{-2}U).
\end{equation}
Write $V$ as an extension
\begin{equation}\label{ext-}
0\longrightarrow M\longrightarrow V\longrightarrow \Lambda M^{-1}\longrightarrow 0.
\end{equation}

Consider $\cS^\Lambda_{\alpha_k^-}$, with $\alpha_k\neq\alpha_m,\alpha_M$ if $d_U$ is even or just $\alpha_k\neq\alpha_M$ if $d_U$ is odd (the case of $\alpha_k=\alpha_m$ and $d_U$ even is easy and will be covered in Proposition \ref{S-alphak01} below).
In these cases, the map $\gamma':M\to M^{-1}U$ as given in \eqref{gamma'} is not an isomorphism. Let $D=\divisor(\gamma')$ be its divisor and consider the structure sheaf $\mathcal{O}_D$ of $X$ restricted to $D$.
Recall that 
one has then a short exact sequence of sheaves
$
0\longrightarrow\mathcal{O}(-D)\longrightarrow\mathcal{O}\stackrel{r(D)}{\longrightarrow}\mathcal{O}_D\longrightarrow 0$ where $r(D)$ is the truncation map.

The following is analogous to Propositions 3.17 and 3.18 of \cite{gothen-oliveira:2011}.
\begin{proposition}\label{prop:s-gamma}
  Let $(V,\gamma)\in\cS^\Lambda_{\alpha_k^-}$. Then there is a well defined section $s_\gamma\in H^0(D,\Lambda^{-1}U)$ given by $s_\gamma=r(D)(\gamma|_M)$.
  Moreover $s_\gamma=0$ if and only if the extension (\ref{ext-}) is trivial, $V\cong M\oplus\Lambda M^{-1}$ and, with respect to this decomposition, 
$$\gamma=\begin{pmatrix}
     \gamma' & 0 \\
     0 & \gamma''
\end{pmatrix},$$
where $\gamma'$ is defined in \eqref{gamma'} and $\gamma''\in H^0(X,\Lambda^{-2}M^2U)\setminus\{0\}$. 
\end{proposition}

We can therefore write $\cS^\Lambda_{\alpha_k^-}$ as a disjoint union
\begin{equation}\label{dijunion}
\cS^\Lambda_{\alpha_k^-}=\cS^{\Lambda,0}_{\alpha_k^-}\sqcup\cS^{\Lambda,1}_{\alpha_k^-}
\end{equation}
where
\begin{itemize}
 \item $\cS^{\Lambda,0}_{\alpha_k^-}$ is the space of pairs in $\cS^\Lambda_{\alpha_k^-}$ with $s_\gamma=0$ ($\Leftrightarrow$ \eqref{ext-} splits);
 \item $\cS^{\Lambda,1}_{\alpha_k^-}$ is the space of pairs in $\cS^\Lambda_{\alpha_k^-}$ with $s_\gamma\neq 0$  ($\Leftrightarrow$ \eqref{ext-} does not split).
\end{itemize}

Consider also the subvariety $\mathcal C(M,\gamma')$ of $H^0(D,\Lambda^{-1}U)\setminus\{0\}\times H^0(X,\Lambda^{-2}U^2)\setminus\{0\}$ whose elements $(q,\eta)$ satisfy the equation $q^2+\eta|_D=0$. 
Then we have a free $\C^*$-action on $\mathcal C(M,\gamma')$ given by $\lambda\cdot(q,\eta)=(\lambda q,\lambda^2\eta)$ and we denote the quotient by
\begin{equation}\label{Q}
\mathcal Q(M,\gamma')=\mathcal C(M,\gamma')/\C^*.
\end{equation}

\begin{proposition}\label{S-alphak01}
Let $\alpha_k\neq\alpha_M$. Let $d$ and $d_U$ be as in Assumption \ref{assump:dU-d>2g-2}. Then in the decomposition \eqref{dijunion} of $\cS^\Lambda_{\alpha_k^-}$:
\begin{enumerate}
\item $\cS^{\Lambda,0}_{\alpha_k^-}$ is projective bundle over $\cS^{\Lambda,0}_{\alpha_k^-}\to\cN_{d-\alpha_k^+}(1,d-\alpha_k)$, with fiber $\mathbb{P}^{d_U-2\alpha_k-g}$.
\item $\cS^{\Lambda,1}_{\alpha_k^-}$ is such that there is a morphism $\cS^{\Lambda,1}_{\alpha_k^-}\to\cN_{d-\alpha_k^+}(1,d-\alpha_k)$ whose fiber over $(M,\gamma')$ is isomorphic to $\mathcal Q(M,\gamma')$ as defined in (\ref{Q}). Moreover, if $d_U$ is even, $\cS^{\Lambda,1}_{\alpha_m^-}=\emptyset$.
\end{enumerate}
\end{proposition}
\proof
Consider first $\alpha_k\neq\alpha_m,\alpha_M$ if $d_U$ is even or just $\alpha_k\neq\alpha_M$ if $d_U$ is odd.
From Proposition \ref{prop:s-gamma}, we know that if $(V,\gamma)\in\cS^{\Lambda,0}_{\alpha_k^-}$, then $V\cong M\oplus\Lambda M^{-1}$, with $\deg(M)=d-\alpha_k$ by \eqref{eq:MdestabS-}, and $\gamma=\begin{pmatrix}
          \gamma' & 0\\
	  0 & \gamma''
         \end{pmatrix}$ with $\gamma'\in H^0(X,M^{-2}U)$ and $\gamma''\in H^0(X,\Lambda^{-2}M^2U)$, both  non-zero.
From this we get a map 
\begin{equation}\label{eq:S0-toN(1,d-alpha)}
\cS^{\Lambda,0}_{\alpha_k^-}\longrightarrow\cN_{d-\alpha_k^+}(1,d-\alpha_k),
\end{equation} sending $(V,\gamma)$ to $(M,\gamma')$. Notice that there is no ambiguity on the choice of the pair $(M,\gamma')$ because $\deg(M)>\deg(\Lambda M^{-1})$. The fiber of this map over $(M,\gamma')$ is given by $\mathbb{P}H^0(X,\Lambda^{-2}M^2U)$, which has
constant dimension, equal to $-2\alpha_k+d_U+1-g$. This completes the description of $\cS^{\Lambda,0}_{\alpha_k^-}$.

The case of $\cS^{\Lambda,1}_{\alpha_k^-}$ is a direct adaptation to the fixed determinant version of \cite[Proposition 3.22]{gothen-oliveira:2011}.

When $\alpha_k=\alpha_m$ and $d_U$ is even, the map $\gamma':M\to M^{-1}U$ is an isomorphism. 
So, in this case, the divisor $D=\divisor(\gamma')$ is not defined. However, Proposition 3.14 and Corollary 3.15 of \cite{gothen-oliveira:2011} readily adapt to our fixed determinant version, showing that $\cS^\Lambda_{\alpha_m^-}$ is a projective bundle over $\cN_{(d_U/2)^-}(1,d_U/2)$ with fiber $\PP^{2(d_U-d)-g}$. So, by (3) of Proposition \ref{prop:Nalpha(1,d)}, it is in fact a disjoint union of $2^{2g}$ copies of $\PP^{2(d_U-d)-g}$. 
Indeed, \cite[Proposition 3.14]{gothen-oliveira:2011} shows that in this case \eqref{ext-} always splits, so we uniformize our statement by saying that $\cS^\Lambda_{\alpha_m^-}=\cS^{\Lambda,0}_{\alpha_m^-}$ and that $\cS^{\Lambda,1}_{\alpha_m^-}$ is empty.
\endproof

From the previous result and from Proposition \ref{prop:Nalpha(1,d)}, we have
\begin{equation}\label{eq:dimS-0}
\dim\cS^{\Lambda,0}_{\alpha_k^-}=2(d_U-d)-g.
\end{equation}
On the other hand, one can prove (cf. \cite[Corollary 3.23]{gothen-oliveira:2011}) that the dimension of 
$\mathcal Q(M,\gamma')$ is independent of $(M,\gamma')$ and equal to $2(d_U-d)-g$, so Proposition \ref{prop:Nalpha(1,d)} tells us again that 
\begin{equation}\label{eq:dimS-1}
\dim\cS^{\Lambda,1}_{\alpha_k^-}=3d_U-4d+2\alpha_k-g.
\end{equation}
Since $\alpha_k\geq\alpha_m\geq d-d_U/2$, it follows that, for every critical value $\alpha_k$, $\dim\cS^{\Lambda,1}_{\alpha_k^-}\geq\dim\cS^{\Lambda,0}_{\alpha_k^-}$, so:
\begin{corollary}\label{dimSalphak-(2,d)}
Let $\alpha_k\neq\alpha_M$. Then, for every critical value $\alpha_k$, each connected component of $\cS^\Lambda_{\alpha_k^-}$ has dimension less or equal than $3d_U-4d+2\alpha_k-g$.
\end{corollary}

With these descriptions, the codimensions of $\cS^\Lambda_{\alpha_k^\pm}$ are easily obtained:

\begin{proposition}\label{Nalpha(d,2) connected} Let $\alpha_k<\alpha_M=d/2$ be a critical value. Then  the codimensions of $\cS^\Lambda_{\alpha_k^\pm}$ in $\cN^\Lambda_{\alpha_k^\pm}$ are such that $\codim\cS^\Lambda_{\alpha_k^+}>3g-3$ and $\codim\cS^\Lambda_{\alpha_k^-}>g-1$.
\end{proposition}

\section{The singular locus}

In this section we describe explicitly the singular locus of the moduli spaces $\cN^\Lambda_\alpha$ and this description will allow us to give a lower estimate of its codimension. We will also conclude that $\cN^{\Lambda,sm}_\alpha$, introduced in Definition \ref{df:smooth}, is indeed the smooth locus of $\cN^\Lambda_\alpha$.

\begin{definition}\label{def:PS alpha}
Let $\alpha\leq d/2$.
\begin{enumerate}
\item Let $\mathcal{PS}_\alpha\subset\cN^\Lambda_\alpha$ be the locus of strictly $\alpha$-polystable quadratic pairs for which the $\alpha$-destabilizing subbundle is of type \textbf{(B)} (cf. Definition \ref{def:types}).
\item Let $\mathcal{NS}_\alpha\subset\cN^\Lambda_\alpha$ be the locus of $\alpha$-stable but non-simple pairs.
\end{enumerate}
\end{definition}

The following lemma will be used below.
\begin{lemma}\label{lemma:bundles in NS and PS not stable}\mbox
If $(V,\gamma)$ is a quadratic pair in $\mathcal{PS}_\alpha$ or in $\mathcal{NS}_\alpha$ then the vector bundle $V$ is not stable. 
\end{lemma}
\proof
If $(V,\gamma)\in\mathcal{PS}_\alpha$, the assertion is clear.
If $(V,\gamma)\in \mathcal{NS}_\alpha$ it is also immediate, because if $V$ was stable, then its only automorphisms were the scalars $\lambda\in\C^*$. The compatibility with $\gamma$ would imply $\lambda=\pm 1$, thus contradicting the non-simplicity of $(V,\gamma)$. 
\endproof

\begin{proposition}\label{prop:codim PS alpha}
For any $\alpha\leq d/2$, and $d$ and $d_U$ as in Assumption \ref{assump:dU-d>2g-2}, $\codim(\mathcal{PS}_\alpha)>4g-5$.
\end{proposition}
\proof
If $(V,\gamma)\in\mathcal{PS}_\alpha$, then by definition, $V\cong L\oplus L'$ with $L'\cong\Lambda L^{-1}$ and $\deg(L)=d/2$, and moreover, $\gamma$ is given, with respect to this decomposition, by
$\gamma=\begin{pmatrix}
          0 & \gamma'\\
	  \gamma' & 0
         \end{pmatrix}$ with $\gamma'\in H^0(X,\Lambda^{-1}U)$ non-zero. These pairs are thus parametrized by $\Jac^{d/2}(X)\times\mathbb{P}H^0(X,\Lambda^{-1}U)$. Since $d_U-d>2g-2$, then $\dim\mathbb{P}H^0(X,\Lambda^{-1}U)=d_U-d-g$, and the result follows from this.
\endproof

Now, recall that the description of the flip loci $\cS^\Lambda_{\alpha_k^-}$, for $\alpha_k<\alpha_M=d/2$, was obtained in a different way than that of $\cS^\Lambda_{\alpha_M^-}$ which we know, from (\ref{eq:SalphaM-}), that it is equal to $\cN_{\alpha_M^-}$. In particular, the spaces $\cS^{\Lambda,0}_{\alpha_k^-}$ were only defined for $\alpha_k<\alpha_M=d/2$.

Let $A$ be the subspace of $\cN^{\Lambda}_{\alpha_M^-}$ given by those quadratic pairs such that the underlying vector bundle is strictly semistable:
\begin{equation}\label{eq:SalphaM-}
A=\{(V,\gamma)\in\cN^{\Lambda}_{\alpha_M^-}\st\exists M\subset V, \deg(M)=d/2\}.
\end{equation}
Obviously, if $d$ is odd, $A$ is empty. 
If $(V,\gamma)\in A$, then there is some $M\subset V$ such that $\deg(M)=d/2$, and the $\alpha_M^-$-stability of $(V,\gamma)$ implies $\gamma(M)\not\subset M^\perp U$.

\begin{remark}\label{rem A}
Notice that this is precisely the analogous situation of \eqref{eq:MdestabS-}, which was related to the description of $\cS^\Lambda_{\alpha_k^-}$, but now we are taking $\alpha_k=\alpha_M=d/2$.
So, if we proceed in the same way as we did above for $\cS^\Lambda_{\alpha_k^-}$ (and if we look at the study of $\cS_{\alpha_k^-}$ carried out in \cite[\S~3.4]{gothen-oliveira:2011}), we immediately conclude that $A$ is described very similarly to $\cS^\Lambda_{\alpha_k^-}$, but by considering  $\alpha_k=\alpha_M$.
\end{remark}

So, we get the following result:

\begin{proposition}\label{prop:A}
Let $A$ be defined as in \eqref{eq:SalphaM-}, and $d$ and $d_U$ as in Assumption \ref{assump:dU-d>2g-2}. Then there is a decomposition 
$A=A^0\sqcup A^1$
where:
\begin{enumerate}
\item $A^0$ fits in the following commutative diagram 
$$\xymatrix{&A^0\ar[r]^(0,35){(V,\gamma)\mapsto (M,\gamma')}\ar[d]_{(V,\gamma)\mapsto (\Lambda M^{-1},\gamma'')}&\cN_{d/2^-}(1,d/2)\ar[d]^{(L,\delta)\mapsto L}\\
&\cN_{d/2^-}(1,d/2)\ar[r]_{(L,\delta)\mapsto\Lambda L^{-1}}&\Jac^{d/2}(X).}$$
So $A^0$ is the fiber product $\cN_{d/2^-}(1,d/2)\times_{\Jac^{d/2}(X)}\cN_{d/2^-}(1,d/2)$ under the given maps.
\item $A^1$ is such that there is a morphism $A^1\to\cN_{d/2^-}(1,d/2)$ whose fiber over $(M,\gamma')$ is isomorphic to $\mathcal Q(M,\gamma')$ as defined in (\ref{Q}).
\end{enumerate}
\end{proposition}

\begin{remark}
It should be remarked that the difference between the descriptions of $A^0$ and of $\cS^{\Lambda,0}_{\alpha_k^-}$ in Proposition \ref{S-alphak01} follows from the fact that for $\alpha_k<\alpha_M$, we always had $\deg(M)>\deg(\Lambda M^{-1})$, so we could always make a consistent choice for the definition of the map \eqref{eq:S0-toN(1,d-alpha)}. In the case of the previous proposition we have $\deg(M)=\deg(\Lambda M^{-1})$, and consequently the description of $A^0$ as a fiber product and not as a projective bundle. 
\end{remark}

Thus, similarly to \eqref{eq:dimS-0} and to \eqref{eq:dimS-1} we have
\begin{equation}\label{eq:A0}
\dim A^0=2(d_U-d)-1 \hspace{.5cm} \text{and} \hspace{.5cm} 
\dim A^1=3(d_U-d)-g.
\end{equation}

\begin{proposition}\label{prop: decomposition of NS alpha}
Let $\alpha<d/2$, and for each critical value $\alpha_k<\alpha_M$, let $\cS^{\Lambda,0}_{\alpha_k^-}$ be defined as in Proposition \ref{S-alphak01}. Let also $A^0$ be as in Proposition \ref{prop:A}. Then: 
$$\mathcal{NS}_\alpha\cong\bigsqcup_{\substack{\alpha_k\text{ critical value }\\ \alpha<\alpha_k<\alpha_M}}\cS^{\Lambda,0}_{\alpha_k^-}\text{  if } d \text{ is  odd}$$ or $$\mathcal{NS}_\alpha\cong\bigsqcup_{\substack{\alpha_k\text{ critical value }\\ \alpha<\alpha_k<\alpha_M}}\cS^{\Lambda,0}_{\alpha_k^-}\sqcup A^0\text{  if } d \text{ is  even}.$$
Furthermore, if $(V,\gamma)\in\mathcal{NS}_\alpha$ then $\Aut(V,\gamma)\cong\Z/2\times\Z/2$.
\end{proposition}

\begin{remark}
This description of $\mathcal{NS}_\alpha$ also holds when $\alpha$ is itself a critical value.
\end{remark}

\proof 
Let $(V,\gamma)\in \mathcal{NS}_\alpha$, and let $f\in\Aut(V,\gamma)\setminus\{\pm 1_V\}$. Then \begin{equation}\label{commute}
\gamma f=((f^t)^{-1}\otimes 1_U)\gamma.
\end{equation}
The coefficients of the characteristic polynomial of $f$ are holomorphic, thus constant, so its eigenvalues are also constant. Let us see that these eigenvalues are different. Indeed, if $f$ has only one eigenvalue $\lambda_0\in\C^*$ then we write $f=\lambda_01_V+f_0$, where $f_0$ is a nilpotent endomorphism of $V$. If $F$ is the line subbundle of $V$ defined as the kernel of $f_0$ then $f_0$ factors through $V/F\cong\Lambda F^{-1}$ and its image lies in $F$. Hence $f_0\in H^0(X,F^2\Lambda^{-1})$. On the other hand, \eqref{commute} implies that $\gamma(F)=0$, the $\alpha$-stability of $(V,\gamma)$ implies $\deg(F)<\alpha$, and $\deg(F^2\Lambda^{-1})<2\alpha-d<0$. The conclusion is that $f_0$ must be zero and $f$ a scalar automorphism of $V$. However, the only scalar automorphisms of $V$ which satisfy \eqref{commute} are precisely $\pm 1_V$, so $f$ cannot be scalar. 
The eigenvalues of $f$ are thus distinct, so we can split $$V\cong F_1\oplus F_2$$ 
as a direct sum of eigenbundles of $f$, such that $F_1F_2\cong\Lambda$ and $F_1\ncong F_2$. 
Suppose that $\deg(F_1)\leq\deg(F_2)$.

Let us now prove that if $(V,\gamma)\in\mathcal{NS}_\alpha$, then $\gamma$ must be generically non-degenerate. Indeed, if $\det(\gamma)=0$, let $N\subset V$ be the line subbundle given by the kernel of $\gamma$. Then (\ref{commute}) shows that either $N\cong F_1$ or $N\cong F_2$. Since we are assuming $\deg(F_1)\leq\deg(F_2)$, we must in fact have by the $\alpha$-stability, $N\cong F_1$ and $\deg(F_1)<\alpha$, and $$\gamma=\begin{pmatrix}
0 & 0\\
0 & \gamma'' 
\end{pmatrix}$$ with $\gamma''\neq 0$. But then $\deg(F_2)>d-\alpha$ and since $\gamma(F_2)\not\subset F_2^\perp U$, this contradits the $\alpha$-stability of $(V,\gamma)$.

So, let $(V,\gamma)\in\mathcal{NS}_\alpha$ with $\gamma$ is genericaly non-degenerate i.e. $\det(\gamma)\neq 0$. In this case, we have that $(\Lambda,\det(\gamma))$ is a $U^2$-quadratic pair of type $(1,d)$. Moreover, $\det(f)$ is an automorphism of $(\Lambda,\det(\gamma))$, so $\det(f)=\pm 1$.
If $\det(f)=1$, (\ref{commute}) shows that,
with respect to the decomposition $V\cong F_1\oplus F_2$,
$$\gamma=\begin{pmatrix}
0 & \gamma'\\
\gamma' & 0 
\end{pmatrix}$$ with $\gamma'\in H^0(X,\Lambda^{-1}U)\setminus\{0\}$. It follows that $\gamma(F_i)\subset F_i^\perp U$ so, by the $\alpha$-semistability of $(V,\gamma)$, we must have $\deg(F_i)=d/2$, hence $(V,\gamma)$ is not $\alpha$-stable, which is a contradiction.
We must then have $\det(f)=-1$. From (\ref{commute}) we conclude that, with respect to the decomposition $V\cong F_1\oplus F_2$,
\begin{equation}\label{eq:nonsimple gamma}
\gamma=\begin{pmatrix}
\gamma' & 0\\
0 & \gamma'' 
\end{pmatrix}
\end{equation} with $\gamma'\in H^0(X,F_1^{-2}U)$ and $\gamma''\in H^0(X,F_1^2\Lambda^{-2}U)$, both non-zero, and that
\begin{equation}\label{eq:non-scalar automorphism}
f=\begin{pmatrix}
 \pm 1 & 0\\
 0 & \mp 1
\end{pmatrix}.
\end{equation}
Recall that we have $F_1\ncong F_2$ (i.e. $F_1$ is not a square root of $\Lambda$). Since $\deg(F_1)\leq\deg(F_2)$, we must have $\deg(F_1)\leq d/2\leq\deg(F_2)$ and $\alpha$-stability implies, $\alpha<\deg(F_1)<d-\alpha$. In other words, there is some critical value $\alpha_k\in (\alpha,d/2]$ such that $\deg(F_1)=\alpha_k\text{ and }\deg(F_2)=d-\alpha_k$.
From \eqref{eq:non-scalar automorphism}, we see that $(V,\gamma)\in\cS^{\Lambda,0}_{\alpha_k^-}$, where $\cS^{\Lambda,0}_{\alpha_k^-}$ is defined in \eqref{dijunion}. 
This gives rise to an isomorphism 
$$\mathcal{NS}_\alpha\cong\bigsqcup_{\substack{\alpha_k\text{ critical value }\\ \alpha<\alpha_k<\alpha_M}}\cS^{\Lambda,0}_{\alpha_k^-}$$ if $d$ is odd or $$\mathcal{NS}_\alpha\cong\bigsqcup_{\substack{\alpha_k\text{ critical value }\\ \alpha<\alpha_k<\alpha_M}}\cS^{\Lambda,0}_{\alpha_k^-}\sqcup A^0$$ if $d$ is even.
From this analysis, it is clear from (\ref{eq:non-scalar automorphism}) that if $(V,\gamma)\in\mathcal{NS}_\alpha$, then $\Aut(V,\gamma)\cong(\Z/2)^2$.
\endproof

Using this description of $\mathcal{NS}_\alpha$ and from \eqref{eq:dimS-0} (as well as from \eqref{eq:A0} if $d$ is even) the following corollary  is immediate.

\begin{corollary}\label{prop:codim NS alpha}
Let $\alpha\leq d/2$ and $d_U-d>2g-2$. Then the locus of $\alpha$-stable but not simple pairs has codimension bigger than $3g-3$.
\end{corollary}

Now we can prove that $\cN^{\Lambda,sm}_\alpha$, introduced in Definition \ref{df:smooth}, is indeed the smooth locus of $\cN^\Lambda_\alpha$. This result will be important in proof of the Torelli type theorem below for these moduli spaces.

\begin{proposition}\label{prop:Nsm is the smooth locus}
For any $\alpha$, $\cN^{\Lambda,sm}_\alpha$ coincides with the smooth locus of $\cN^\Lambda_\alpha$.
\end{proposition}
\proof
Consider first that $\alpha$ is a generic value of the parameter. In this case, $\cN^\Lambda_{\alpha}-\cN^{\Lambda,sm}_{\alpha}=\mathcal{PS}_\alpha\sqcup \mathcal{NS}_\alpha$, so we have to show that the points in $\mathcal{PS}_\alpha$ and in $ \mathcal{NS}_\alpha$ are singular points.

Let $(V,\gamma)\in\mathcal{PS}_\alpha$. Then $V=L\oplus L'$ with $L'\cong\Lambda L^{-1}$ and $\deg(L)=d/2$, and
$\gamma=\begin{pmatrix}
          0 & \gamma'\\
	  \gamma' & 0
         \end{pmatrix}$.

From (\ref{eq:les}), and recalling that $d_U-d>2g-2$ implies that $\mathbb{H}^2(X,C^\bullet(V,\gamma))=0$, it follows that $\dim \mathbb{H}^1(X,C^\bullet(V,\gamma))=3(d_U-d)+\dim\mathbb{H}^0(X,C^\bullet(V,\gamma))$, so the dimension of the infinitesimal deformation space of $(V,\gamma)$ is $\dim \mathbb{H}^1(X,C^\bullet(V,\gamma))/\C=3(d_U-d)+\dim\mathbb{H}^0(X,C^\bullet(V,\gamma))-1$.
Now, the automorphisms of elements in $\cN^{\Lambda,sm}_\alpha$ are only in $\Z/2$, hence we always have that the corresponding infinitesimal automorphism space $\mathbb{H}^0$ is zero. However, for $(V,\gamma)$ as above, for each $\lambda\in\C$, 
$\psi=\begin{pmatrix}
          \lambda & 0\\
	  0 & -\lambda
         \end{pmatrix}$ is an element of $\mathbb{H}^0(X,C^\bullet(V,\gamma))$, thus $\dim\mathbb{H}^0(X,C^\bullet(V,\gamma))\geq 1$ (this corresponds to the fact that for the given pair $(V,\gamma)$ the automorphisms of the form $\begin{pmatrix}
          \lambda & 0\\
	  0 & 1/\lambda
         \end{pmatrix}$, $\lambda\in\C^*$ are also allowed). So,
  $$\dim \mathbb{H}^1(X,C^\bullet(V,\gamma))/\C\geq 3(d_U-d)>\dim\cN^\Lambda_\alpha,$$ and thus $(V,\gamma)$ represents a singularity of $\cN^\Lambda_\alpha$.
  
Consider now a pair $(V,\gamma)\in\mathcal{NS}_\alpha$. From Proposition \ref{prop: decomposition of NS alpha}, $\Aut(V,\gamma)\cong(\Z/2)^2$. Here, $\mathbb{H}^0(X,C^\bullet(V,\gamma))=0$, but the presence of finitely many more automorphisms than $\Z/2$ yields an orbifold type of singularity of the moduli space. 

If $\alpha=\alpha_k$ is a critical value, then $\cN_{\alpha_k}-\cN_{\alpha_k}^{sm}=\mathcal{PS}_{\alpha_k}\sqcup \mathcal{NS}_{\alpha_k}\sqcup(\cS^\Lambda_{\alpha_k^+}\cup\cS^\Lambda_{\alpha_k^-})$. This follows from the fact that the elements of $\cS^\Lambda_{\alpha_k^+}$ and $\cS^\Lambda_{\alpha_k^-}$ are strictly $\alpha_k$-polystable with destabilizing subbundle of type $\textbf{(A)}$ and $\textbf{(C)}$ while the ones of $\mathcal{PS}_{\alpha_k}$ are strictly $\alpha_k$-polystable with destabilizing subbundle of type $\textbf{(B)}$. So, by \cite[Proposition 2.17]{gothen-oliveira:2011}, $\mathcal{PS}_\alpha$ must be disjoint of $\cS^\Lambda_{\alpha_k^+}\cup\cS^\Lambda_{\alpha_k^-}$. The elements of $\mathcal{NS}_{\alpha_k}$ are $\alpha_k$-stable, so this case is clear. Thus we just have to show that if $(V,\gamma)\in\cS^\Lambda_{\alpha_k^+}\cup\cS^\Lambda_{\alpha_k^-}$ then it represents a singular point. But such $(V,\gamma)$ is strictly $\alpha_k$-polystable of type \textbf{(A)} or \textbf{(C)}, and one sees that it represents a singular point just as we did above for type \textbf{(B)}:
for type \textbf{(A)} we have that the $1$-dimensional family of endomorphisms of type 
$\psi=\begin{pmatrix}
          \lambda & 0\\
	  0 & 0
         \end{pmatrix}$ ($\lambda\in\C$) belong to $\mathbb{H}^0(X,C^\bullet(V,\gamma))$, and for type \textbf{(C)} the same occurs with the $1$-dimensional family 
$\psi=\begin{pmatrix}
          0 & 0\\
	  0 & \lambda
         \end{pmatrix}$, $\lambda\in\C$.
\endproof

For a critical value $\alpha_k<\alpha_M$, let
\begin{equation}\label{eq:flip intersect smooth}
\cS^{\Lambda,sm}_{\alpha_k^\pm}=\cN^{\Lambda,sm}_{\alpha_k^\pm}\cap\cS^\Lambda_{\alpha_k^\pm}.
\end{equation}

Proposition \ref{Nalpha(d,2) connected} tells us that the codimensions of $\cS^{\Lambda,sm}_{\alpha_k^\pm}$ in $\cN^{\Lambda,sm}_{\alpha_k^\pm}$ satisfy the following:
\begin{equation}\label{eq:codimflipinsmooth}
\codim\cS^{\Lambda,sm}_{\alpha_k^-}>g-1\hspace{.5cm}\text{and}\hspace{.5cm}\codim\cS^{\Lambda,sm}_{\alpha_k^+}>3g-3.
\end{equation}
\begin{proposition}\label{N--S-=N+-S+}
If $\alpha_k$ is a critical value, then $\cN^{\Lambda,sm}_{\alpha_k^-}-\cS^{\Lambda,sm}_{\alpha_k^-}
  \cong\cN^{\Lambda,sm}_{\alpha_k}\cong\cN^{\Lambda,sm}_{\alpha_k^+}-\cS^{\Lambda,sm}_{\alpha_k^+}$.
\end{proposition}
\proof
First, notice that the subvarieties $\mathcal{PS}_\alpha$ are independent of $\alpha$, whence $(V,\gamma)\notin\mathcal{PS}_{\alpha_k^\pm}$ if and only if $(V,\gamma)\notin\mathcal{PS}_{\alpha_k}$.
On the other hand, it is also clear, from the description of $\mathcal{NS}_\alpha$ given in Proposition \ref{prop: decomposition of NS alpha}, that if $(V,\gamma)\notin\mathcal{NS}_{\alpha_k^\pm}$ if and only if $(V,\gamma)\notin\mathcal{NS}_{\alpha_k}$.
The result follows from these two facts and from the definition of the flip loci.
\endproof

The next result shows that the singular locus of $\cN^\Lambda_{\alpha}$ has, in general, high codimension.

\begin{proposition}\label{prop:codimension of singular locus} Let $\alpha<d/2$, $g\geq 2$ and let $d$ and $d_U$ be as in Assumption \ref{assump:dU-d>2g-2}. If $\alpha$ is a generic value, then $\cN^\Lambda_{\alpha}$ is smooth outside of a subset of codimension larger than $3g-3$. If $\alpha=\alpha_k$ is a critical value, then $\cN^\Lambda_{\alpha_k}$ is smooth outside of a subset of codimension larger than $g-1$.
\end{proposition}
\proof
For a generical value $\alpha$, we have that $\cN^\Lambda_{\alpha}-\cN^{\Lambda,sm}_{\alpha}=\mathcal{PS}_\alpha\sqcup \mathcal{NS}_\alpha$, so $\codim(\cN_{\alpha}-\cN_{\alpha}^{sm})$ is at least the minimum of $\codim(\mathcal{PS}_\alpha)$ and $\codim(\mathcal{NS}_\alpha)$, so Proposition \ref{prop:codim PS alpha} and Corollary \ref{prop:codim NS alpha} prove (1).

For a critical value $\alpha_k$, notice that $\cN^\Lambda_{\alpha_k}-\cN^{\Lambda,sm}_{\alpha_k}=\mathcal{PS}_{\alpha_k}\sqcup\mathcal{NS}_{\alpha_k}\sqcup(\cS^\Lambda_{\alpha_k^+}\cup\cS^\Lambda_{\alpha_k^-})$. So (2) follows by Proposition \ref{prop:codim PS alpha}, Corollary \ref{prop:codim NS alpha} and also Proposition \ref{Nalpha(d,2) connected}.
\endproof

\begin{lemma}\label{lemma:codim of A}
Let $A$ be the locus defined in \eqref{eq:SalphaM-} and let $A^{sm}=A\cap\cN^{\Lambda,sm}_{\alpha_M^-}$. Then $\codim(A^{sm})\geq g-1$.
\end{lemma}
\proof
Taking into account Remark \ref{rem A}, the stated codimension is given just by replacing $\alpha_k$ in $\codim(\cS^{\Lambda,sm}_{\alpha_k^-})$ in \eqref{eq:codimflipinsmooth} by $\alpha_M=d/2$.
\endproof

\begin{remark}\label{rem:smoothnessoflargecv}
Notice that if $d$ is odd, then $\mathcal{PS}_\alpha$ is empty for every $\alpha$ and the same holds for $\mathcal{NS}_{\alpha_M^-}$ (the case of $\mathcal{NS}_{\alpha_M^-}$ follows from Proposition \ref{prop: decomposition of NS alpha} and from the fact that $A$ is also empty). So, in this case, $\mathcal{N}^\Lambda_{\alpha_M^-}$ is indeed smooth. However, we cannot conclude that $\cN^\Lambda_\alpha$ is smooth for $\alpha$ such that there is a critical value $\alpha_k$ with $\alpha<\alpha_k<\alpha_M^-$, because $\mathcal{NS}_\alpha$ is non-empty.
\end{remark}

\section{Topological properties}

Let $\cM^{\Lambda,s}=\cM^s(2,\Lambda)$ be the stable locus of the moduli space $\cM^\Lambda$ of rank $2$ polystable vector bundles with fixed determinant $\Lambda$. Of course, if $d$ is odd, then $\cM^{\Lambda,s}=\cM^\Lambda$.

\begin{proposition}\label{prop:projective bundle last moduli}
Let $\pi:\cN^{\Lambda}_{\alpha_M^-}\to\cM^\Lambda$ be the map defined in Proposition \ref{fiber d/2}. 
Then:
\begin{enumerate}
\item $\pi(\cN^{\Lambda}_{\alpha_M^-}-\cN^{\Lambda,sm}_{\alpha_M^-})\subset\cM^\Lambda-\cM^{\Lambda,s}$;
\item Let $A$ be the subspace of $\cN^{\Lambda}_{\alpha_M^-}$ given by those quadratic pairs such that the underlying vector bundle is strictly semistable, defined in \eqref{eq:SalphaM-}, and let $A^{sm}=A\cap\cN^{\Lambda,sm}_{\alpha_M^-}$. Let $p$ be the restriction of $\pi$ to $\cN^{\Lambda,sm}_{\alpha_M^-}-A^{sm}$. Then $p:\cN^{\Lambda,sm}_{\alpha_M^-}-A^{sm}\to\cM^{\Lambda,s}$ is a projective bundle, with fiber isomorphic to $\mathbb P^{3(d_U-d-g)+2}$.
\end{enumerate}   
\end{proposition}
\proof
Recall that $\cN^\Lambda_{\alpha_M^-}-\cN^{\Lambda,sm}_{\alpha_M^-}=\mathcal{PS}_{\alpha_M}\sqcup \mathcal{NS}_{\alpha_M}$ where $\mathcal{PS}_\alpha$ and $\mathcal{NS}_\alpha$ are defined in Propositions \ref{prop:codim PS alpha} and \ref{prop:codim NS alpha}, respectively. From Lemma \ref{lemma:bundles in NS and PS not stable} we know that if $V$ occurs in a pair $(V,\gamma)$ in $\mathcal{PS}_{\alpha_M^-}$ or in $\mathcal{NS}_{\alpha_M^-}$, then it is not stable. Furthermore, if $(V,\gamma)\in\mathcal{NS}_{\alpha_M^-}$ then the $\alpha_M^-$-stability ensures that $V$ is polystable (cf. Proposition \ref{fiber d/2}). Hence $\pi(V,\gamma)=V\in\cM^\Lambda-\cM^{\Lambda,s}$.
This proves (1). Obviously, if $d$ is odd, this item carries no information (cf. Remark \ref{rem:smoothnessoflargecv}). 

It is clear that $p(\cN^{\Lambda,sm}_{\alpha_M^-}-A)\subset\cM^{\Lambda,s}$ and, if $V\in\cM^{\Lambda,s}$,  from (1) we know that  $\mathbb{P}H^0(S^2V^*\otimes U)=\pi^{-1}(V)\in\cN^{\Lambda,sm}_{\alpha_M^-}-A$. Then, together with Proposition \ref{fiber d/2}, this yields (2).
\endproof

\begin{proposition}\label{prop:Nsm dense in N}
For any $\alpha$, and $d,d_U$ as in Assumption \ref{assump:dU-d>2g-2}, $\cN^{\Lambda,sm}_\alpha$ is dense in $\cN^\Lambda_\alpha$.
\end{proposition}
\proof
Fix an arbitrary point in $\cN^\Lambda_\alpha-\cN^{\Lambda,sm}_\alpha$ represented by a pair $(V,\gamma)$. We will prove the proposition by considering separately the cases where $\alpha$ is generic or critical.
Suppose $\alpha$ is generic. Then $(V,\gamma)$ belongs to $\mathcal{PS}_\alpha\sqcup \mathcal{NS}_\alpha$, and Lemma \ref{lemma:bundles in NS and PS not stable} states that for every $(V,\gamma)\in\cN^\Lambda_\alpha-\cN^{\Lambda,sm}_\alpha$, $V$ is not stable.
Now, from \cite[Proposition 2.6]{narasimhan-ramanan:1975} that we can deform $V$ into a stable vector bundle $V'$ without changing the determinant. More precisely, there is a non-singular and irreducible scheme $T$ parametrizing a family $\mathcal{V}=\{V_t\}_{t\in T}$ of rank $2$ and degree $d$ vector bundles over $X$ such that $\bigwedge\nolimits^2V_t=\Lambda$, for all $t$, $V_{t_0}=V$ for some $t_0\in T$ and the subspace $\{t\in T\st V_t\text{ stable}\}$ is non-empty and dense in $T$. In other words, $\mathcal{V}\to T\times X$ is a vector bundle such that $\mathcal{V}|_{\{t\}\times X}$ satisfies the above conditions. 
 
Let $p_1:T\times X\to T$ and $p_2:T\times X\to X$ be the projections, and consider the direct image sheaf $R^0p_{1*}(S^2\mathcal{V}^*\otimes p_2^*U)$. Since $d_U-d>2g-2$, the fibers of $R^0p_{1*}(S^2\mathcal{V}^*\otimes p_2^*U)$ have constant (positive) dimension, so it is a locally free sheaf, thus a vector bundle over $T$.
We conclude therefore that the corresponding projective bundle $\mathbb{P}R^0p_{1*}(S^2\mathcal{V}^*\otimes p_2^*U)$ is a connected base for a family of $U$-quadratic pairs which contains the given pair $(V,\gamma)$ and for which there is a non-empty dense subspace consisting of pairs such that the underlying vector bundle is stable. The first paragraph of the proof implies then that on this non-empty dense subspace, all quadratic pairs are $\alpha$-stable and simple. From the universal property of the coarse moduli space  $\cN_\alpha$ (cf. \cite[Theorem 1.6]{gomez-sols:2000}), there is a morphism $\mathbb{P}R^0p_{1*}(S^2\mathcal{V}^*\otimes p_2^*U)\to\cN_\alpha$ which factors through $\cN^\Lambda_\alpha\subset\cN_\alpha$ and this provides a deformation of $(V,\gamma)$ to an $\alpha$-stable and simple pair i.e. a pair in $\cN^{\Lambda,sm}_\alpha$.  This proves the proposition for generic $\alpha$.

If $\alpha=\alpha_k$ is a critical value, then we have two obvious continuous maps 
$\pi_{\pm}:\cN_{\alpha_k^{\pm}}^\Lambda\to\cN_{\alpha_k}^\Lambda$. 
From the definition of the flip loci $\cN_{\alpha_k}^\Lambda=\pi_-(\cN_{\alpha_k^-}^\Lambda)\cup\pi_+(\cN_{\alpha_k^+}^\Lambda)$ 
and $\pi_\pm(\cN_{\alpha_k^\pm}^\Lambda)=\cN_{\alpha_k}^\Lambda-\cS_{\alpha_k^\mp}^\Lambda$.
From above, we see that $\pi_\pm(\cN^{\Lambda,sm}_{\alpha_k^\pm})$ is dense in $\pi_\pm(\cN^{\Lambda}_{\alpha_k^\pm})$, thus $\pi_-(\cN^{\Lambda,sm}_{\alpha_k^-})\cup\pi_+(\cN^{\Lambda,sm}_{\alpha_k^+})$ is dense in $\cN_{\alpha_k}^\Lambda$. On the other hand, from Proposition \ref{N--S-=N+-S+}, $\cN^{\Lambda,sm}_{\alpha_k}=\pi_-(\cN^{\Lambda,sm}_{\alpha_k^-}-\cS^{\Lambda,sm}_{\alpha_k^-})=\pi_+(\cN^{\Lambda,sm}_{\alpha_k^+}-\cS^{\Lambda,sm}_{\alpha_k^+})$, and \eqref{eq:codimflipinsmooth} tells us, since the spaces are smooth, that $\cN^{\Lambda,sm}_{\alpha_k^\pm}-\cS^{\Lambda,sm}_{\alpha_k^\pm}$ is dense in $\cN^{\Lambda,sm}_{\alpha_k^\pm}$. So we have $\cN_{\alpha_k}^\Lambda=\pi_-(\cN_{\alpha_k^-}^\Lambda)\cup\pi_+(\cN_{\alpha_k^+}^\Lambda)=\overline{\pi_-(\cN_{\alpha_k^-}^{\Lambda,sm})}\cup\overline{\pi_+(\cN_{\alpha_k^+}^{\Lambda,sm})}=\overline{\cN_{\alpha_k}^{\Lambda,sm}}$, as claimed.
\endproof

From \cite[Theorem 3.5]{gomez-sols:2000} we know that under some conditions $\cN^\Lambda_\alpha$ is irreducible. The previous proposition gives us a different proof of the same result.

\begin{corollary}\label{irreducible}
Let $g\geq 2$ and $d,d_U$ satisfying Assumption \ref{assump:dU-d>2g-2}. Then $\cN^\Lambda_\alpha$ is irreducible.
\end{corollary}
\proof
Since the moduli space $\cM^\Lambda$ of rank two stable bundles with fixed determinant is irreducible \cite[Remark 5.9]{newstead:2012}, it follows from (2) of Proposition \ref{prop:projective bundle last moduli} that $\cN_{\alpha_M^-}^{\Lambda,sm}-A^{sm}$ is irreducible. Since $\cN_{\alpha_M^-}^{\Lambda,sm}$ is smooth and $g\geq 2$, then from Lemma \ref{lemma:codim of A} we conclude that $\cN_{\alpha_M^-}^{\Lambda,sm}$ is irreducible as well. The same reasoning, now using  Proposition \ref{N--S-=N+-S+} and \eqref{eq:codimflipinsmooth}, shows that $\cN^{\Lambda,sm}_\alpha$ is irreducible for any $\alpha$. Finally the previous proposition completes the proof.
\endproof

Our techniques can also be used to compute the fundamental group of the smooth locus of $\cN^\Lambda_\alpha$.

\begin{corollary}\label{pi1=0}
Let $g\geq 2$ and $d,d_U$ as in Assumption \ref{assump:dU-d>2g-2}. Then $\cN^{\Lambda,sm}_\alpha$ is simply-connected.
\end{corollary}
\proof
The moduli space $\cM^{\Lambda,s}$ of rank two stable bundles with fixed determinant is simply-connected (cf. \cite[Corollary 2]{newstead:1967}), so from (2) of Proposition \ref{prop:projective bundle last moduli} we conclude that $\pi_1(\cN^{\Lambda,sm}_{\alpha_M^-}-A)=0$. As $\codim(A)$ is positive and $\cN^{\Lambda,sm}_{\alpha_M^-}$ is smooth, we have $\pi_1(\cN^{\Lambda,sm}_{\alpha_M^-})=0$.
Now, suppose that $\pi_1(\cN^{\Lambda,sm}_{\alpha_k^+})=0$ for some $\alpha_k<\alpha_M$.
Since $\cN^{\Lambda,sm}_{\alpha_k^+}$  is smooth, using again Proposition \ref{N--S-=N+-S+} and \eqref{eq:codimflipinsmooth}, we have
$$0=\pi_1(\cN^{\Lambda,sm}_{\alpha_k^+})=\pi_1(\cN^{\Lambda,sm}_{\alpha_k^+}-\cS^{\Lambda,sm}_{\alpha_k^+})=\pi_1(\cN^{\Lambda,sm}_{\alpha_k})=\pi_1(\cN^{\Lambda,sm}_{\alpha_k^-}-\cS^{\Lambda,sm}_{\alpha_k^-})=\pi_1(\cN^{\Lambda,sm}_{\alpha_k^-}),$$ whence the result.
\endproof

Notice that one cannot, a priori, make any statement about the fundamental group of $\cN^\Lambda_\alpha$, because the presence of singularities may have great influence on its the topology.

\section{Mixed Hodge structures and the Torelli Theorem}

\subsection{Hodge structures}

Here we recall the basics on the theory of (mixed) Hodge structures on complex algebraic varieties.
The main references are the papers \cite{deligne:1970,deligne:1971,deligne:1974} by Deligne and also the books \cite{peters-steenbrink:2008,voisin:2002}.

\begin{definition}
A \emph{(pure) Hodge structure of weight $k$} is a pair $(V_\Z,V_\C)$ where $V_\Z$ is a free abelian group and $V_\C=V_\Z\otimes\C$ is its complexification, such that there is a decomposition into subspaces $V_\C=\bigoplus_{p+q=k}V^{p,q}$ with $\overline{V^{q,p}}\cong V^{p,q}$. There is a corresponding Hodge filtration given by $F^p(V_\C)=\bigoplus_{r\geq p}V^{r,s}$.
A \emph{polarized Hodge structure} is a Hodge structure $(V_\Z,V_\C)$ together with a bilinear nondegenerate map $\theta:V_\Z\otimes V_\Z\to\Z$, such that $V^{p,q}$ and $V^{p',q'}$ are orthogonal under $\theta_\C:V_\C\otimes V_\C\to\C$, unless $p'=k-p$ and $q'=k-q$.
A \emph{morphism between two Hodge structures} $V_\C$, $V'_\C$ is a linear map $f:V_\C\to V'_\C$ compatible with the Hodge filtration, and a morphism between two polarized Hodge structures $V_\C$ and $V'_\C$ is a morphism of Hodge structures preserving the polarizations.
\end{definition}

For example, given a smooth projective variety $Z$ (hence compact Kähler, with Kähler form $\omega$) of dimension $n$, its cohomology has a pure Hodge structure given by $H^i(Z,\C)=\bigoplus_{p+q=i}H^{p,q}(Z)$. This Hodge structure is polarized, with the polarization given by the cup product:
\begin{equation}\label{eq:polarization of smooth projective variety}
\theta(\alpha,\beta)=\int_Z\alpha\cup\beta\cup\omega^{n-k}, 
\end{equation}
for $\alpha,\beta\in H^k(Z,\Z)$.

It is well-known that to any Hodge structure of weight $1$ (in fact of any odd weight) there is associated a complex torus. To see this, consider the inclusion $V_\Z\hookrightarrow V_\C$ and compose this with the projection $p_{0,1}$ of $V_\C$ onto, say, $V^{0,1}$. The quotient $V^{0,1}/p_{0,1}(V_\Z)$ is a complex torus. Conversely, from a complex torus $V/\Gamma$ we obtain a Hodge structure of weight $1$ by taking $V_\C=V\oplus\overline V$ and $V_\Z=\Gamma$ embedded in $V_\C$. This procedure applied to $H^1(X,\C)$, where $X$ is our smooth projective curve, gives rise to the Jacobian of $X$: $H^{0,1}(X)/H^1(X,\Z)\cong\Jac(X)$.
If, in this construction, we consider the polarization $\theta$, then we obtain a $(1,1)$-form in the complex torus, which is integer and nondegenerate, hence Kähler. This is in fact the first Chern class of an ample line bundle on the torus. In other words, we get a polarized abelian variety. When applied to the Jacobian of $X$ as above, we obtain the canonical polarization in $\Jac(X)$, such that the divisor of the corresponding ample line bundle over $\Jac(X)$ is the theta divisor, $\Theta$.
So, having the polarized Hodge structure $(H^1(X),\theta)$, with $\theta$ given by \eqref{eq:polarization of smooth projective variety}, is equivalent to having the polarized Jacobian $(\Jac(X),\Theta)$, which gives therefore an alternative way to state the classical Torelli theorem:

\begin{theorem}[\textbf{Torelli theorem}]\label{thm:Torelli}
If $X$ and $X'$ are two smooth projective curves such that $(H^1(X),\theta)$ and $(H^1(X'),\theta')$ are isomorphic as polarized Hodge structures, then $X\cong X'$.
\end{theorem}

As we saw above, the cohomology of a smooth projective variety has a pure Hodge structure. This is not true anymore for varieties which are not smooth or projective. In these cases we have to consider mixed Hodge structures, which we now define.

\begin{definition}
Let $V_\mathbb{Z}$ be a free abelian group. A \emph{mixed Hodge structure} over $V_\mathbb{Z}$ consists of an ascending weight filtration $W$ on the rational finite dimensional vector space $V_\mathbb{Q}=V_\mathbb{Z}\otimes \mathbb{Q}$ and a descending Hodge filtration $F$ on $V_\C$ such that $F$ induces a pure Hodge filtration of weight $k$ on each rational vector space $\mathrm{Gr}_k^W V=V_k/V_{k-1}$. 
\end{definition}

Define $V^{p,q}=\mathrm{Gr}_F^p(\mathrm{Gr}^W_{p+q}V)_\C$. A mixed Hodge structure is pure of weight $k$ if $\mathrm{Gr}_l^W V=0$ for $l\neq k$.
If $Z$ is a complex algebraic variety (not necessarily projective or smooth), Deligne has shown in \cite{deligne:1970,deligne:1971,deligne:1974} that the integral cohomology $H^k(Z)$ and the cohomology with compact support $H_c^k(Z)$ both carry natural (mixed) Hodge structures. If $Z$ is a smooth projective variety then this is the pure Hodge structure previously mentioned.
The notion of isomorphism of (mixed) Hodge structures and isomorphism of polarized (mixed) Hodge structures are the obvious ones.

Now, we will consider the moduli spaces $\cN^{\Lambda,sm}_{\alpha}$ which are smooth, quasi-projective, but not projective. So, a priori, their cohomology has a mixed Hodge structure. However, has we shall see below, at least the first cohomology groups (the ones we shall use), have pure Hodge structure, so we will not really use in practice the notion of mixed Hodge structure. 

\subsection{Cohomology groups of $\cN^{\Lambda,sm}_{\alpha}$}

The following lemmas are proved in \cite[Lemma 6.1.1]{arapura-sastry:2000} and \cite[Corollary 6.1.2]{arapura-sastry:2000}, respectively.
 
\begin{lemma}\label{lemma:arapura-sastry} Let $M$ be a smooth variety and let $Z\subset M$ be a closed subvariety of codimension $k$. Then
$H^i(M,\Z)\cong H^i(M-Z,\Z)$ for all $i<2k-1$.
\end{lemma}

\begin{lemma}\label{lemma:arapura-sastry II} Let $M$ be a projective variety and let $U\subset M$ be a smooth Zariski open subset, with $\codim(M-U)=k$. Then
$H^i(U,\Z)$ has a pure Hodge structure of weight $i$, for every $i<k-1$.
\end{lemma}

With these results we can now easily compute the first integral cohomology groups of $\cN^{\Lambda,sm}_{\alpha}$. Write $$H^*(\cN^{\Lambda,sm}_\alpha)=H^*(\cN^{\Lambda,sm}_\alpha,\mathbb{Z}).$$
Recall that $\cM^{\Lambda,s}$ denotes the stable locus of the moduli space $\cM^\Lambda$. 

\begin{proposition}\label{Prop:cohomology groups}
Under Assumption \ref{assump:dU-d>2g-2}, let $g\geq 4$ and $\alpha<d/2$. Then, for every $\alpha<\alpha_M$, we have, $H^1(\cN^{\Lambda,sm}_\alpha)=0$, $H^2(\cN^{\Lambda,sm}_\alpha)\cong\Z\oplus\Z$ and $H^3(\cN^{\Lambda,sm}_\alpha)\cong H^1(X)$.
\end{proposition}
\proof
It is known \cite{newstead:1967} that  $H^1(\cM^{\Lambda,s})=0$, $H^2(\cM^{\Lambda,s})\cong \Z$ and $H^3(\cM^{\Lambda,s})\cong H^1(X)$.
So (2) of Proposition \ref{prop:projective bundle last moduli} implies that
$H^1(\cN^{\Lambda,sm}_{\alpha_M^-}-A)=0$, $H^2(\cN^{\Lambda,sm}_{\alpha_M^-}-A)\cong \Z\oplus\Z$ and $H^3(\cN^{\Lambda,sm}_{\alpha_M^-}-A)\cong H^1(X)$.
Since $g\geq 4$, the result for $\alpha=\alpha_M^-$  follows from Lemmas \ref{lemma:codim of A} and \ref{lemma:arapura-sastry} and for any $\alpha$ from Proposition \ref{N--S-=N+-S+}, \eqref{eq:codimflipinsmooth} and Lemma \ref{lemma:arapura-sastry}. Of course, the result for $H^1$ also follows from Corollary \ref{pi1=0}.
\endproof

\begin{corollary}
If $g\geq 4$ and $\alpha<d/2$, then the Picard group of $\cN^{\Lambda,sm}_\alpha$ is isomorphic to $\Z\oplus\Z$.
\end{corollary}
\proof
Consider first $\cN^{\Lambda,sm}_{\alpha_M^-}$, and identify $\Pic(\cN^{\Lambda,sm}_{\alpha_M^-})$ with $H^1(\cN^{\Lambda,sm}_{\alpha_M^-},\mathcal{O}^*)$. From the previous proposition we have $H^1(\cN^{\Lambda,sm}_{\alpha_M^-})=0$ and  $H^2(\cN^{\Lambda,sm}_{\alpha_M^-})\cong\Z\oplus\Z$, thus the exponential sequence yields the exact sequence
\begin{equation}\label{eq:exactseq for Picard}
0\longrightarrow H^1(\cN^{\Lambda,sm}_{\alpha_M^-},\mathcal{O})\longrightarrow\Pic(\cN^{\Lambda,sm}_{\alpha_M^-})\longrightarrow\Z\oplus\Z.
\end{equation}
As $H^1(\cN^{\Lambda,sm}_{\alpha_M^-})=0$, then $H^1(\cN^{\Lambda,sm}_{\alpha_M^-},\C)=0$. Moreover, from Proposition \ref{prop:codimension of singular locus}, $\codim(\cN_{\alpha_M^-}-\cN^{\Lambda,sm}_{\alpha_M^-})\geq 3g-3$ so, since $g\geq 4$, Lemma \ref{lemma:arapura-sastry II} says that $H^1(\cN^{\Lambda,sm}_{\alpha_M^-})$ is pure, and we have $H^1(\cN^{\Lambda,sm}_{\alpha_M^-},\C)=H^{1,0}(\cN^{\Lambda,sm}_{\alpha_M^-})\oplus H^{0,1}(\cN^{\Lambda,sm}_{\alpha_M^-})$, with $\overline{H^{0,1}}(\cN^{\Lambda,sm}_{\alpha_M^-})\cong H^{1,0}(\cN^{\Lambda,sm}_{\alpha_M^-})$. It follows that
 $H^{0,1}(\cN^{\Lambda,sm}_{\alpha_M^-})=0$ i.e. $H^1(\cN^{\Lambda,sm}_{\alpha_M^-},\mathcal{O})=0$, therefore the map $\Pic(\cN^{\Lambda,sm}_{\alpha_M^-})\to\Z\oplus\Z$ in (\ref{eq:exactseq for Picard}) is injective.
In order to see that it is also surjective, note that, since $H^2(\cN^{\Lambda,sm}_{\alpha_M^-})=\Z\oplus\Z$, we have $b_2=2$, $b_2$ denoting the second Betti number of $\cN^{\Lambda,sm}_{\alpha_M^-}$. Using again Lemma \ref{lemma:arapura-sastry II}, we see that 
$H^2(\cN^{\Lambda,sm}_{\alpha_M^-},\C)=H^{2,0}(\cN^{\Lambda,sm}_{\alpha_M^-})\oplus H^{1,1}(\cN^{\Lambda,sm}_{\alpha_M^-})\oplus H^{0,2}(\cN^{\Lambda,sm}_{\alpha_M^-})$, with $H^{0,2}(\cN^{\Lambda,sm}_{\alpha_M^-})\cong\overline{H^{2,0}}(\cN^{\Lambda,sm}_{\alpha_M^-})$, so 
\begin{equation}\label{eq:secondBetti}
2=b_2=h^{1,1}+2h^{2,0}
\end{equation} where $h^{i,j}=\dim H^{i,j}(\cN^{\Lambda,sm}_{\alpha_M^-})$ are the Hodge numbers. Since $\cN$ is projective it has a (positive) line bundle whose Chern class is a $(1,1)$-class. Restricting this line bundle to $\cN^{\Lambda,sm}_{\alpha_M^-}$ shows then that $h^{1,1}\geq 1$, so by (\ref{eq:secondBetti}) we must have  $h^{1,1}=2$. In other words, all classes in $H^2(\cN^{\Lambda,sm}_{\alpha_M^-})=\Z\oplus\Z$ come from $(1,1)$-classes, i.e. $\Pic(\cN^{\Lambda,sm}_{\alpha_M^-})\to\Z\oplus\Z$ is surjective.

Hence $\Pic(\cN^{\Lambda,sm}_{\alpha_M^-})\cong\Z\oplus\Z$, and by Proposition \ref{N--S-=N+-S+} and \eqref{eq:codimflipinsmooth}, we conclude that $\Pic(\cN^{\Lambda,sm}_\alpha)\cong\Z\oplus\Z$ for every $\alpha< d/2$.
\endproof

\subsection{The Torelli Theorem for $\cN^\Lambda_\alpha$}
In the next two sections, we shall consider families of curves, so our moduli spaces $\cN^{\Lambda}_\alpha$ will depend on which curve we are considering. Hence, for a curve $X$, we will write $\cN^{\Lambda}_{X,\alpha}$ to emphasize that we are considering quadratic pairs over $X$.

\begin{proposition}
Under Assumption \ref{assump:dU-d>2g-2}, let $g\geq 5$. Then $H^3(\cN^{\Lambda,sm}_{X,\alpha})$ is naturally polarized, and the isomorphism $H^3(\cN^{\Lambda,sm}_{X,\alpha})\cong H^1(X)$ respects the polarizations.
\end{proposition}
\proof This is similar to \cite{arapura-sastry:2000} and \cite{munoz:2009}.
Let $H_1$ and $H_2$ be two positive generators of $\Pic(\cN^{\Lambda,sm}_{X,\alpha})\cong\Z\oplus\Z$. For every $\alpha$, $\cN^\Lambda_{X,\alpha}$ is projective, thus naturally polarized  i.e. there are $a,b\in\Z$ such that $H=aH_1+bH_2$ is a polarization of $\cN^\Lambda_{X,\alpha}$. If $k=3(d_U-d)-1$ is the dimension of $\cN^\Lambda_{X,\alpha}$, take a generic $(k-3)$-fold hyperplane intersection $Z\subset\cN^{\Lambda,sm}_{X,\alpha}$. Hence $Z$ is smooth, $\dim Z=3$ and, by Proposition \ref{prop:codimension of singular locus}, we see that $Z$ is in fact projective.
Now, let
\begin{equation}\label{eq:polarization}
  \hat\theta_X:H^3(\cN^{\Lambda,sm}_{X,\alpha})\otimes H^3(\cN^{\Lambda,sm}_{X,\alpha})\longrightarrow\Z,\hspace{.5cm}\hat\theta(\alpha\otimes\beta)=\langle\alpha\cup\beta,[Z]\rangle=\int_Z \alpha\cup\beta.
\end{equation}
We prove that this is a polarization as follows. Take a generic $(k-4)$-fold hyperplane section $W\subset\cN^{\Lambda,sm}_{X,\alpha}$. Proposition \ref{prop:codimension of singular locus}, and the fact that $g\geq 5$, tells us again that $W$ is projective and smooth. Then there is a version of the Lefschetz hyperplane theorem compatible to our case (cf. \cite[Theorem 6.1.1]{arapura-sastry:2000}) which when applied to $\cN^\Lambda_{X,\alpha}$, $\cN^{\Lambda,sm}_{X,\alpha}$ and $W$ yields $H^3(W)\cong H^3(\cN^{\Lambda,sm}_{X,\alpha})$ and, by the Hard Lefschetz theorem, we have an isomorphism
$H^3(W)\xrightarrow{\cong}H^5(W)\cong H^3(W)^*$. This map coincides with $\hat\theta_X$, proving thus that it is non-degenerate.

Let now $\pi:\chi\to T$ be a family of generic curves of genus $g$, parametrized by an irreducible base variety $T$, such that the natural map $T\to \mathfrak{M}_g$ to the moduli space $\mathfrak{M}_g$ of genus $g$ curves is dominant. Let $q:\mathcal J^d\to T$ be the universal Jacobian. Then we consider the universal moduli space
$p:\widetilde\cN^T_\alpha(2,d)\to\mathcal J^d$ such that $p^{-1}(X,\Lambda)=\cN^{\Lambda,sm}_{X,\alpha}$.
Let $f:\widetilde\cM^{s,T}(2,d)\to\mathcal J^d$ be the universal moduli space of stable vector bundles, such that $f^{-1}(X,\Lambda)=\cM^{\Lambda,s}_X$.
There is then a forgetful map 
$\widetilde\cN^T_{\alpha_M^-}(2,d)\to\widetilde\cM^{s,T}(2,d)$, and pulling back the relative ample generator of 
$\widetilde\cM^{s,T}\to\mathcal J^d$, we obtain an element $H_2\in\Pic(\cN^{\Lambda,sm}_{X,\alpha_M^-})$, defined in the family $\pi$. On the other hand, since $\widetilde\cN^T_{\alpha_M^-}(2,d)\to\widetilde\cM^{s,T}(2,d)$ is, by Proposition \ref{prop:projective bundle last moduli}, outside of a subset of codimension bigger than $2$, a projective bundle, there is another element $H_1\in\Pic(\cN^{\Lambda,sm}_{X,\alpha_M^-})$, well-defined in families.
The construction of the flips can also be done in families (cf. \cite{gothen-oliveira:2011}) so, for any $\alpha$, $\Pic(\cN^{\Lambda,sm}_{X,\alpha})\cong\Z\oplus\Z$ with the generators defined in families.
The construction of the polarization $\hat\theta_X$ in \eqref{eq:polarization} goes then through in families as well, hence defines a polarization in $R^3p_*\underline\Z$, which we denote by  $\hat\theta$.

Now, let $\theta$ be the standard polarization on $R^1\pi_*\underline\Z$ given by the cup product - see \eqref{eq:polarization of smooth projective variety}. Then, we know that there is an isomorphism $R^3p_*\underline\Z\cong q^*R^1\pi_*\underline\Z$
which, together with the polarization $\hat\theta$, yields a polarization $\theta'$ in $q^*R^1\pi_*\underline\Z$. Since the map $T\to \mathfrak{M}_g$ is dominant, then Lemma 8.1.1 of \cite{arapura-sastry:2000} asserts that $\theta'$ is a positive integer multiple of $\theta$ hence so is $\hat\theta$.
Finally, restricting to $X$ again, we conclude that $\hat\theta_X$ is a positive multiple of $\theta_X$, hence it determines a unique primitive polarization on $H^3(\cN^{\Lambda,sm}_{X,\alpha})$. So the isomorphism $H^3(\cN^{\Lambda,sm}_{X,\alpha})\cong H^1(X)$ respects the polarizations.
\endproof

\begin{theorem}\label{torelli for quad pairs}
Let $X$ and $X'$ be smooth projective curves of genus $g,g'\geq 5$, $\Lambda$ and $\Lambda'$ line bundles of degree $d$ and $d'$ on $X$ and $X'$, respectively, and $U$ and $U'$ line bundles of degree $d_U$ and $d_U'$ on $X$ and $X'$, respectively, such that $d,d_U$ and $d',d_U'$ are as in Assumption \ref{assump:dU-d>2g-2}. Let $\cN^\Lambda_{X,\alpha}$  be the moduli space of $\alpha$-polystable $U$-quadratic pairs on $X$, and define similarly $\cN^{\Lambda'}_{X',\alpha}$.
If $\cN^\Lambda_{X,\alpha}\cong \cN^{\Lambda'}_{X',\alpha}$ then $X\cong X'$. 
\end{theorem}
\proof
Due to Proposition \ref{prop:Nsm is the smooth locus}, the isomorphism $\cN^\Lambda_{X,\alpha}\cong \cN^{\Lambda'}_{X',\alpha}$ induces an isomorphism $\cN^{\Lambda,sm}_{X,\alpha}\cong \cN^{\Lambda',sm}_{X',\alpha}$ and therefore an isomorphism between the corresponding $H^3$. This isomorphism preserves the corresponding polarizations $\hat\theta_X$ and $\hat\theta_{X'}$ as constructed in the previous theorem, due to their definition. So, we have
$(H^3(\cN^{\Lambda,sm}_{X,\alpha}),\hat\theta_X)\cong (H^3(\cN^{\Lambda',sm}_{X',\alpha}),\hat\theta_{X'})$. The previous theorem says then that $(H^1(X),\theta_X)\cong (H^1(X'),\theta_{X'})$ where $\theta_X$ and $\theta_{X'}$ are the corresponding standard polarizations given by \eqref{eq:polarization of smooth projective variety}. The usual Torelli Theorem \ref{thm:Torelli} completes the proof.\endproof

\subsection{The non-fixed determinant case}

Recall that $\cN^\Lambda_\alpha:=\mathrm{det}^{-1}(\Lambda)$ where $\det$ is the map given by (\ref{eq:detmap}).
Since the fibers of $\det$ are all isomorphic, we have that $\det:\cN_{X,\alpha}(2,d)\to\Jac^d(X)$ is a fibration, where, as before, we write $\cN_{X,\alpha}(2,d)$ to emphasize that we are considering the moduli space of rank $2$, degree $d$, quadratic pairs over $X$.
We will extend the results of the previous sections to the moduli space $\cN_{X,\alpha}(2,d)$ of $\alpha$-polystable $U$-quadratic pairs over $X$, using the techniques of \cite{munoz:2009}.

The following lemma, proved in \cite[Lemma 7.1]{munoz:2009}, will be needed in the proof of the Torelli theorem for $\cN_{X,\alpha}(2,d)$.

\begin{lemma}
Let $M$ be a projective connected variety and $f:M\to Y$ a map to another quasi-projective variety such that $f^*:H^k(Y)\to H^k(M)$ is zero for all $k>0$. Then $f$ is constant.
\end{lemma}

\begin{theorem}\label{torelli for quad pairs non fixed det} Let $X$ and $X'$ be smooth projective curves of genus $g,g'\geq 5$ and let $U$ and $U'$ be line bundles of degree $d_U$ and $d_U'$ on $X$ and $X'$, respectively. Let $d$ and $d'$ be two integers numbers, such that $d,d_U$ and $d',d_U'$ are as in Assumption \ref{assump:dU-d>2g-2}. If $\cN_{X,\alpha}(2,d)\cong\cN_{X',\alpha}(2,d')$, then $X\cong X'$. 
\end{theorem}
\proof
Let $F:\cN_{X,\alpha}(2,d)\stackrel{\cong}{\longrightarrow}\cN_{X',\alpha}(2,d')$ be the given isomorphism. Let $\cN^{sm}_{X,\alpha}(2,d)$ be the subspace defined  precisely in the same way as we did in Definition \ref{df:smooth} for the fixed determinant case. Define similarly $\cN^{sm}_{X',\alpha}(2,d')$. The same arguments as in Proposition \ref{prop:Nsm is the smooth locus} show that $\cN^{sm}_{X,\alpha}(2,d)$ is the smooth locus of $\cN_{X,\alpha}(2,d)$. So $F$ restricts to an isomorphism (still denoted by $F$) of the corresponding smooth locus
$F:\cN^{sm}_{X,\alpha}(2,d)\stackrel{\cong}{\to}\cN^{sm}_{X',\alpha}(2,d')$.
Fix $\Lambda_0\in\Jac^d(X)$ and take the composition $$f:\cN_{X,\alpha}^{\Lambda_0,sm}\hookrightarrow\cN^{sm}_{X,\alpha}(2,d)\stackrel{F}{\longrightarrow}\cN^{sm}_{X',\alpha}(2,d')\longrightarrow\Jac^{d'}(X').$$
From Proposition \ref{Prop:cohomology groups} we know that $H^1(\cN^{\Lambda_0,sm}_{X,\alpha})=0$, so $f^*:H^1(\Jac^{d'}(X'))\to H^1(\cN^{\Lambda_0,sm}_{X,\alpha})$ is the zero map. Moreover, the cohomology of $\Jac^{d'}(X')$ is generated by $H^1(\Jac^{d'}(X'))$, so $f^*:H^k(\Jac^{d'}(X'))\to H^k(\cN^{\Lambda_0,sm}_{X,\alpha})$ is the zero map for every $k>0$. The previous lemma implies that $f$ is constant, therefore $F(\cN^{\Lambda_0,sm}_{X,\alpha})\subset\cN^{\Lambda_0',sm}_{X,\alpha}$, for some $\Lambda_0'\in\Jac^{d'}(X')$. Since $\cN^{\Lambda_0,sm}_{X,\alpha}$ and $\cN^{\Lambda_0',sm}_{X,\alpha}$ are dense in $\cN^{\Lambda_0}_{X,\alpha}$ and in $\cN^{\Lambda_0'}_{X,\alpha}$, respectively, we have $F(\cN^{\Lambda_0}_{X,\alpha})\subset\cN^{\Lambda_0'}_{X,\alpha}$. In the same way, we see that $F^{-1}(\cN^{\Lambda_0'}_{X,\alpha})\subset\cN^{\Lambda_0}_{X,\alpha}$. In other words, we have an isomorphism $$F|_{\cN^{\Lambda_0}_{X,\alpha}}:\cN^{\Lambda_0}_{X,\alpha}\stackrel{\cong}{\longrightarrow}\cN^{\Lambda_0'}_{X,\alpha}$$
and Theorem \ref{torelli for quad pairs} yields the result.
\endproof

\bibliographystyle{amsplain}

\end{document}